\documentclass[twocolumn]{autart}
\usepackage{graphicx}      
\usepackage{natbib}

\usepackage{amsmath,amssymb,mathrsfs,dsfont}

\newtheorem{assumption}{Assumption}

\usepackage{graphicx}
\graphicspath{{./fig/}}
\usepackage[usenames,dvipsnames]{color}
\usepackage{epstopdf}
\usepackage{subfig}
\usepackage{siunitx}
\usepackage[americanresistors,americaninductors]{circuitikz}

\usetikzlibrary{calc}
\usepackage{tikz}
\usetikzlibrary{shapes,arrows}
\usetikzlibrary{decorations.pathmorphing}
\ctikzset{bipoles/thickness=0.7}
\ctikzset{bipoles/length=0.8cm}
\ctikzset{bipoles/diode/height=.375}
\ctikzset{bipoles/diode/width=.3}
\ctikzset{tripoles/thyristor/height=.8}
\ctikzset{tripoles/thyristor/width=1}
\ctikzset{bipoles/vsourceam/height/.initial=.7}
\ctikzset{bipoles/vsourceam/width/.initial=.7}
\tikzstyle{every node}=[font=\small]
\tikzstyle{every path}=[line width=0.7pt,line cap=round,line join=round]

\tikzstyle{block} = [draw, fill=blue!10, rectangle, minimum height=4em, minimum width=8em]
\tikzstyle{smallblock} = [draw, fill=blue!10, rectangle, minimum height=2em, minimum width=2em]
\tikzstyle{sum} = [draw, fill=blue!10, circle, node distance=1cm]
\tikzstyle{input} = [coordinate]
\tikzstyle{output} = [coordinate]
\tikzstyle{pinstyle} = [pin edge={to-,thin,black}]

\newcommand\oprocendsymbol{\hbox{$\square$}}
\newcommand\oprocend{\relax\ifmmode\else\unskip\hfill\fi\oprocendsymbol}

\newcommand*\real[0]{\mathbb R}       

\newcommand*\mycircle[0]{\mathbb S}
\renewcommand{\vec}[1]{\mathbf{#1}}

\usepackage{mathtools}
\usepackage{enumitem}
\usepackage{xcolor}
\usepackage[bbgreekl]{mathbbol}
\usepackage{amsfonts}
\DeclareMathAlphabet{\mathpzc}{OT1}{pzc}{m}{it}

\DeclareSymbolFontAlphabet{\mathbb}{AMSb}
\DeclareSymbolFontAlphabet{\mathbbl}{bbold}

\providecommand{\norm}[1]{\lVert#1\rVert}

\allowdisplaybreaks

\begin{document}
\begin{frontmatter}

\title{Grid-forming Control for Power Converters based on Matching of Synchronous Machines \thanksref{footnoteinfo}} 

\thanks[footnoteinfo]{The first two authors contributed equally.\\ This work was partially funded by the European Union's Horizon 2020 research and innovation programme under grant agreement N$^\circ$ 691800, ETH Z\"urich funds, and the SNF Assistant Professor Energy Grant \#160573. This article reflects only the authors' views and the European Commission is not responsible for any use that may be made of the information it contains. A preliminary version of part the results in this paper has been presented at the IFAC-Workshop on Distributed Estimation and Control in Networked Systems, September 8-9, 2016 Tokyo, Japan \citep{TJ-CA-FD:16}.}
\author[First]{Catalin Arghir}\ead{carghir@control.ee.ethz.ch} ~
\author[First]{Taouba Jouini}\ead{tjouini@control.ee.ethz.ch} ~
\author[First]{Florian D\"orfler}\ead{dorfler@control.ee.ethz.ch}
 
\address[First]{Automatic Control Laboratory at the Swiss Federal Institute of Technology (ETH) Z\"urich, Switzerland.}

\begin{abstract} 
We consider the problem of grid-forming control of power converters in low-inertia power systems. 
Starting from an average-switch three-phase {power converter} model, we draw parallels to a synchronous machine (SM) model {and} propose a novel converter control strategy which dwells upon the main characteristic of a SM: the presence of an internal rotating magnetic field. 
In particular, we augment the converter system with a virtual oscillator whose frequency is driven by {the} DC-side voltage measurement and which sets the converter pulse-width-modulation signal, {thereby} achieving exact matching between the converter in closed-loop and the SM dynamics. 
We then provide a sufficient condition asserting existence, uniqueness, and global asymptotic stability of a shifted equilibrium, all in a rotating coordinate frame attached to the virtual oscillator angle. 
By actuating the DC-side input of the converter we are able to enforce this condition and provide additional inertia and damping. 
In this framework, we illustrate strict incremental passivity, droop, and power-sharing properties which are compatible with conventional power system operation requirements.
We subsequently adopt disturbance-decoupling and droop techniques to design additional control loops that regulate the DC-side voltage, as well as AC-side frequency and amplitude, while in the end evaluating them with numerical experiments.

\end{abstract}


\end{frontmatter}

\section{Introduction}
\label{sec: intro}

The electrical power system is currently undergoing significant changes in its structure and mode of operation due to a major shift in generation technology from synchronous machines (SMs) to  power electronics-based DC/AC converters, or simply {\em inverters}. As opposed to SMs, which store kinetic energy in their rotor moment of inertia, these devices are on the one hand designed with little or no built-in energy storage capacity, while on the other hand actuated at much faster time scales. 
SMs, with their large rotational inertia, self-synchronizing physics and associated controls, act as safeguards against faults and disturbances -- all of which are absent in {\em low-inertia systems} with a dominant share of distributed and variable renewable sources interfaced through inverters. Hence, the proper control of inverters is regarded as one of the key challenges when massively integrating renewable energy sources \citep{BK-BJ-YZ-VG-PD-BMH-BH:17,JT-SVD-DSC:16,GD-TP-PP-XK-FC-XG:15}.

Converter control strategies are classified in two groups. While there is no universally accepted definition, inverters are usually termed {\em grid-following} if their controls are designed for a stiff grid, and they deliver power at the stiff AC grid frequency usually measured through a phase-locked loop (PLL). Otherwise, these converters are termed {\em grid-forming} when they are assigned to interact with a non-stiff grid similarly as SMs do by balancing kinetic and electrical energy in such a way that a frequency consensus is achieved. 
A low-inertia system cannot be operated with only grid-following units. With this in mind, we review the literature on grid-forming control.

The inherent self-synchronizing property of SMs has inspired controllers such as {\em droop} and {\em virtual synchronous machines} (VSMs) \citep{torres2013virtual,VK-SDH-HKZ:11,MW-SH-PV-KV:09,SDA-SJA:13,YC-RH-DT-HPB:11
,QCZ-GW:11}. These controllers are designed to emulate the behavior of SM models of different degree of fidelity and are based on measurements of AC quantities such as injected power, frequency, and amplitude. For example, {\em inverse droop} and related VSM control strategies measure the AC frequency through a PLL and accordingly adapt the converter power injection based on a simple SM swing equation model. The latter is encoded in a micro-controller whose outputs are tracked by the converter modulation signal typically through a cascaded control architecture. For these, and other VSM implementations, the time delays resulting from measuring and processing AC quantities render control often ineffective \citep{ENTSOE2016,HB-TI-YM:14,GD-TP-PP-XK-FC-XG:15}.
Droop control can also be implemented by measuring the injected power and by adapting accordingly the converter frequency \citep{JMG-MC-TLL-PLC:13c}, but its applicability is limited to inductive grids and with a possibly narrow region of attraction \citep{sinha2017uncovering,FD-JWSP-FB:14a,de2018bregman}. Additionally, the inverter's DC-side storage element is often not included in the model, nor in the control design, which, in our view, misses a key insight: namely, that the DC bus voltage can reflect power imbalance and serve as valuable feedback signal.
Finally, alternative control strategies {employ} {\em nonlinear virtual oscillators} fed by AC current measurements \citep{BBJ-SVD-AOH-PTK:14a,sinha2017uncovering,colombino2017global}. For these strategies global stability certificates are known, but their design and analysis is quite involved (as a result, no controllers for regulation of amplitudes and frequency are known thus far) and their compatibility with SMs is unclear to this date. 

Another set of literature relevant to our methodology is {\em passivity-based control} (PBC) \citep{AJvdS:96} and {\em interconnection and damping assignment} (IDA) \citep{ortega2004interconnection}. Their application to DC/DC converters \citep{escobar1999hamiltonian,zonetti2014globally}, AC/DC converters \citep{perez2004passivity}, and power systems in general \citep{YS-PT:14,FZO+13} suggests a physically insightful analysis based on shaping the energy and dissipation functions. 
As we will further see, our analysis relies also on a characterization of the power system steady-state specifications \citep{arghir2016steady,gross2017steady} which restrict the class of admissible controllers.

Our main {contributions are} three-fold. 
First, we propose a novel grid-forming control strategy that matches the electromechanical energy exchange pattern in SMs. This is achieved by augmenting the converter dynamics with an internal model of a harmonic oscillator whose frequency tracks the value of the DC-side voltage measurement. This voltage-driven oscillator is then assigned to drive the converter's pulse-width-modulation cycle, thereby assuring that the closed-loop converter dynamics exactly {\em match} the SM dynamics, whereas the DC voltage serves as the key control and imbalance signal akin to the SM's angular velocity. 
Based on a Lyapunov approach we provide a sufficient condition certifying existence, uniqueness, and global asymptotic stability of driven equilibria, in a coordinate frame attached to the virtual oscillator angle. By actuating the DC-side input current we are able to satisfy this condition. We also preserve strict incremental passivity, droop, and power-sharing properties of the closed-loop system. Our approach is grounded in foundational control methods, while being systematically extensible to PBC and IDA designs. Additionally, the key DC voltage signal is readily available while all other approaches rely on extensive processing of AC measurements.
Second, building on the proposed matching controller, we further design overarching control loops that regulate the DC voltage, AC frequency, and AC amplitude. This is done by pursuing an approach based on disturbance decoupling, which performs asymptotic output voltage amplitude tracking, while rejecting the load current seen as a measurable disturbance. We then suggest extensions based on employing PBC and voltage-power droop control strategies, which have been previously investigated in various settings. 
Third and finally, we evaluate the performance and robustness of our designs by comparing them in numerical experiments of single and multi-converter scenarios.

The remainder of the paper is organized as follows. Section \ref{sec: models and control} introduces the models and the control objectives. Section~\ref{sec: matching control} proposes the matching controller and derives its properties. Section \ref{sec: outer loops} designs the regulation and disturbance-decoupling controllers. Section \ref{sec: case study} presents a numerical case study, and Section \ref{sec: conclusions} concludes the paper.

\section{The Three-Phase Converter Model, Synchronous Machine Model, \& their Analogies} 
\label{sec: models and control}

\subsection{Preliminaries {and coordinate transformations}}
\label{subsec: preliminaries}

In this paper $\boldsymbol{I} = \left[\begin{smallmatrix} 1 & 0 \\ 0 & 1 \end{smallmatrix} \right]$ {denotes} the identity and $\boldsymbol{J} = \left[\begin{smallmatrix} 0 & -1 \\ 1 & 0 \end{smallmatrix} \right]$ denotes the rotation by $\tfrac{\pi}{2}$ in $\real^2$, {while} $\mathbf{e}_2=\begin{bsmallmatrix} 0 \\ 1 \end{bsmallmatrix}$ is a natural basis vector in $\real^2$. We denote by $\norm{\,\cdot\,}$ the standard Euclidean norm for vectors or the induced norm for matrices.

 {The three-phase AC system is assumed to be symmetrical, namely all passive elements} have equal values for each phase element. Due to this symmetry, any three-phase quantity $z_{abc}\in\real^3$ is assumed to satisfy $\begin{bsmallmatrix} 1 & 1 & 1 \end{bsmallmatrix}z_{abc} = 0$; see Remark~\ref{remark: zero sequence}. We consider a coordinate transformation to distinguish between the component along the span of the vector $\begin{bsmallmatrix} 1 & 1 & 1 \end{bsmallmatrix}^\top\in \mathbb{R}^3$, which we denote by $z_\gamma\in\real$ and the other two components {$z_{\alpha\beta}\in\real^2$ lying on the associated orthogonal complement called the $\alpha\beta$-frame:}
\begin{equation}
\label{eq: Talphabeta}
\begin{bmatrix}{z}_{\alpha\beta}\\z_\gamma\end{bmatrix} = \begin{smallmatrix}\sqrt{{2}/{3}}\end{smallmatrix} \begin{bsmallmatrix}1 & -\frac{1}{2} & -\frac{1}{2} \\ 0 & {\frac{\sqrt3}{2}} & -{\frac{\sqrt3}{2}} \\ \frac{1}{\sqrt{2}} & \frac{1}{\sqrt{2}} & \frac{1}{\sqrt{2}} \end{bsmallmatrix}z_{abc} \,.
\end{equation}
Given a reduced three-phase quantity $z_{\alpha\beta}$ and an angle $\theta \in \mathbb S^{1}$, we define the {$dq$-coordinate transformation $(z_{\alpha\beta},\theta)\mapsto z_{dq}\in\real^2$, via $\boldsymbol{R}_{\theta} = \begin{bsmallmatrix}
\cos{\theta} & -\sin{\theta} 
\\
\sin{\theta} & \cos{\theta}
\end{bsmallmatrix}$, as} 
\begin{equation}
\label{$dq$-transf}
z_{dq}= \boldsymbol{R}_{\theta}^\top z_{\alpha\beta} \,.
\end{equation}%
{Consequently, we have that a sinusoidal steady-state solution of the form $\dot z_{\alpha\beta}^\star = \omega^*{\boldsymbol{J}}z_{\alpha\beta}^\star$, with associated frequency $\omega^{*}$, is mapped to an equilibrium $\dot z^{*}_{dq}=0$ in the $dq$-frame whose transformation angle satisfies $\dot\theta^\star = \omega^{*}$.} 
%
Throughout this article, a variable denoted $z_{dq}^*$ or $z_{\alpha\beta}^\star$ is used to represent a steady-state solution induced by {exogenous} inputs, e.g., {load parameters or} set-points.

\subsection{Three-Phase DC/AC Converter Model}
\label{subsec: modeling}

We {start by reviewing} the standard average-switch model of a {three-phase, two-level, voltage source inverter in $\alpha\beta$-coordinates. See \citep{yazdani2010voltage} for a comprehensive study. The model is described by a continuous-time system whose main feature is the nonlinearity captured by the modulation (switching) block, as depicted in Figure \ref{fig:circuit diagram} below.

\begin{figure}[h!]
\begin{minipage}{\columnwidth}
\centering
\begin{center}
\resizebox{8.5cm}{3cm}{
\begin{circuitikz}[american voltages]
\draw
(0,0) to [short, *-] (3,0)
(5,3) to [american current source, i>=$i_{l}$] (5,0) 
(3.9,3) to [short, *-] (5,3)
(3.9,0) to [short, *-] (5,0)
(3,0) to (3.9,0)
(3,3) to (3.9,3)
(0.7,0) to [open, v^<=${v}_x$] (0.7,3) 
(0,3) 
to [short,*-, i=$i_{\alpha\beta}$] (0.7,3) 
(0.6,3) to [R, l=$R$] (1.6,3) 
to [L, l=$L$] (2.6,3) 
(2.6,3) to (2.8,3)
to [short,*-] (2.8,2) 
(2.8,1) to [C, l=$C$] (2.8,2) 
(2.8,1) to [short] (2.8,0)
(3.9,0) to [R, l_=$G$] (3.9,3) 
(2.7,3) to (3,3)
(2.95,3) to [open, v^>=${v}_{\alpha\beta}$] (2.95,0); 
\draw
(-2,-0.5) to (-2,3.5)
(-0.5,1.5) node[pigbt] (pigbt){} 
(-2,3.5) to (0,3.5)
(0,3.5) to (0,-0.5)
(-2,-0.5) to (0,-0.5);
\draw
(-2,0) to [short, *-] (-4.3,0)
(-4.3,0) to [short, *-] (-4,0)
(-5.5,0) to (-4.3,0)
(-5.5,3) to (-4.3,3)
(-5.5,0) to [american current source, i>=$i_{dc}$] (-5.5,3) 
(-4.3,3) to [R,l=$G_{dc}$] (-4.3,0) 
(-3.1,3) to [C, l=$C_{dc}$] (-3.1,0) 
(-2,3) to [short, *-, i<=$i_{x}$] (-3,3)
(-4.3,3) to [short, *-] (-3,3)
(-3.1,3) to [short, *-] (-3,3)
(-5.2,3) to [open, v^>=${v}_{dc}$] (-5.2,0);
\end{circuitikz}
}
\end{center}
\caption{Circuit diagram of a $3$-phase DC/AC converter} 
\label{fig:circuit diagram}
\end{minipage}
\end{figure}
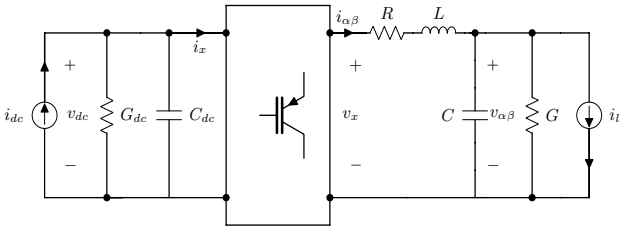

The DC circuit consists of a {controllable} current source {$i_{dc}\in\real$} in parallel with a {capacitance} $C_{dc}>0$ and a conductance $G_{dc}>0$. The DC-side switching current is denoted by $i_{x}\in\real$, while $v_{dc}\in\real$ represents the voltage across the DC capacitance.
The AC circuit contains at each phase an inductance $L>0$ in series with a resistance $R>0$ 
connected to a shunt capacitance $C>0$ {and shunt conductance $G>0$}. Here $v_{\alpha\beta}\in\real^2$ denotes the AC voltage across the output capacitor. 

The dissipative elements $G_{dc}$, $G$, and $R$ model the parasitic losses in the converter. Furthermore, $i_{\alpha\beta}\in\real^2$ denotes the AC current in the inductors and $v_{x}\in\real^2$ the average AC voltage at the switching node. 
{The inverter model is terminated at its AC-ports with a load current $i_{l}$ drawn from a weak AC grid, which will be made more specific in Assumption~\ref{ass: iL-model}.}

The switching block {is defined as} the average-switch\footnote{For the time scales of interest, we assume a sufficiently high switching frequency that allows us to discard the PWM carrier harmonics and use continuous-time dynamics.} model of a $6$-switch $2$-level inverter {with an associated} complementary pulse-width-modulation (PWM) carrier and a modulation signal  
$m_{\alpha\beta}\in \{x\in\real^2: \norm{x}\leq1 \}\,.$
To preserve energy conservation, the switching block is assumed to be lossless, i.e. it {satisfies} the identities
\[
i_{x} = \frac{1}{2}m_{\alpha\beta}^{\top}i_{\alpha\beta} \\ \,\, \, , \, \,\,
v_{x} =\frac{1}{2}m_{\alpha\beta} v_{dc}\,.
\]

By putting it all together, the inverter model can be written as the following bilinear system
\begin{subequations}%
\begin{align}%
C_{dc}\dot v_{dc} &= -G_{dc} v_{dc}+i_{dc}-\frac{1}{2}m_{\alpha\beta}^{\top} {i_{\alpha\beta}}
\label{eq: DC voltage}
\\
\,
L \dot{i}_{\alpha\beta} &= -R i_{\alpha\beta} - v_{\alpha\beta} + \frac{1}{2}m_{\alpha\beta}{v_{dc}}
\label{eq: inductance current}
\\
\,
C \dot v_{\alpha\beta} &= -G v_{\alpha\beta} + i_{\alpha\beta} -i_{l}
\,.
\end{align}%
\label{eq: inverter dynamics}%
\end{subequations}%
\begin{rem}[Zero sequence]\label{remark: zero sequence}
We will construct the three-phase modulation signal $m_{abc}$ in such a way that $m_{\gamma}=0$ which implies that $v_{x,\gamma}=0$. For a balanced load, it also holds that $i_{l,\gamma}=0$. We are left with the following dynamics for the $\gamma$-subsystem:
\begin{subequations}
\begin{align}
L {\dot {i}_{\gamma}} &= -R i_{\gamma} - v_{\gamma}
\\
C \dot v_{\gamma} &=  -G v_{\gamma} + i_{\gamma}
\,.
\end{align}	%
\label{eq: system dynamics}%
\end{subequations}%
Since \eqref{eq: system dynamics} is an asymptotically stable linear system, the omission of the $\gamma$-component is well justified.
\oprocend
\end{rem}

\subsection{Control objectives}

In this section we map out the control objectives to be achieved via the two main actuation inputs, the modulation signal $m_{\alpha\beta}$ and the DC-side current injection $i_{dc}$.  Broadly speaking we require the following: 

{\em (i) Grid-forming}: 
The objective of {\em grid-forming control} is best defined by mimicking the electromechanical interaction of a SM with the grid rather than prescribing the converter's frequency to track the grid frequency, e.g., via a PLL. The synchronization properties of SMs rely on a particular kinetic to electrical energy exchange pattern. This can be induced in the DC/AC converter by exactly {matching} the SM's dynamics.

{\em (ii) Voltage and amplitude regulation:}  
We intend to exactly regulate $v_{{dc}}$ and $\norm{v_{\alpha\beta}}$ to prescribed references, possibly requiring knowledge of system parameters and full state measurements.
If the load current measurements are uncertain or unknown, we aim instead to achieve a linear {\em droop} characteristic between the converter modulation frequency and its power output.
Such a local droop behavior is known to guarantee power sharing and compatibility with other droop-like controllers in a power system \citep{sinha2017uncovering,FD-JWSP-FB:14a}. 

{\em (iii) Strict incremental passivity}: 
We aim to preserve strict incremental passivity \citep{AJvdS:96} with respect to the AC and DC ports, $u = (i_{dc},-i_{l})$ and $y=(v_{dc}, v_{dq})$, and {relative} to a desired steady-state solution $x^{*} = (v_{dc}^{*},i_{dq}^{*},v_{dq}^{*})$. More precisely, we seek a positive definite storage function that is decreasing along system trajectories, where the system remains strictly incrementally passive after implementing the controller.

In the sequel, we further specify these objectives, in more suitable coordinates, and also consider alternative objectives such as voltage amplitude droop.

\subsection{The Synchronous Machine Model}
\label{subsec: emulation control}

In what follows, we consider a SM model which lends itself useful in designing the matching controller. 
We consider a single-pole-pair, non-salient rotor SM under constant excitation, defined in $\alpha\beta$-frame as in \citep{YS-PT:14}, together with a capacitor at its AC terminal, and described by the state-space model
\begin{subequations}%
\label{eq: sync-gen}
\begin{align}
\dot \theta &= \omega
\\
M \dot\omega&=-D \omega+\tau_m+L_m i_f\begin{bmatrix}
-\sin{\theta}\\ \cos{\theta}
\end{bmatrix}^{\top} {i_{\alpha\beta}}
\label{eq: sync-gen-w}
\\
L_s {\dot{i}_{\alpha\beta}}&=-{R_s} i_{\alpha\beta}- v_{\alpha\beta} - L_m i_f
\begin{bmatrix}
-\sin{\theta}\\ \cos{\theta}
\end{bmatrix} {\omega }
\label{eq: sync-gen-i}
\\
C \dot v_{\alpha\beta}&= -Gv_{\alpha\beta}+i_{\alpha\beta}-i_{l}\, .
\end{align}%
\end{subequations}%
Here, $M>0$ and $D>0$ are the rotor inertia and damping coefficients, $\tau_m$ is the driving mechanical torque, $L_m>0$ is the stator-to-rotor mutual inductance, $L_s>0$ the stator inductance. We denote the rotor angle by $\theta\in\mycircle^1$, its angular velocity by $ \omega \in \real$, the current in the stator winding by $i_{\alpha\beta}\in\real^2$, and the stator resistance by ${R_s}>0$. At its terminals the SM is interfaced to the grid through a shunt capacitor with capacitance $C>0$ and capacitor voltage $v_{\alpha\beta}\in\real^2$, a constant load conductance $G>0$, and the load current extraction denoted by $i_{l}\in\real^{2}$.
The strength of the rotating magnetic field inside the SM is given by the rotor current $i_f$ which is assumed to be regulated to a constant (here, negative) value, as in \citep{YS-PT:14,arghir2016steady}.

Observe the similarities between the inverter model \eqref{eq: inverter dynamics} and the SM model \eqref{eq: sync-gen}. The DC capacitor is analogous to the rotor moment of inertia, while the electrical torque and the electromotive force (EMF) (the rightern-most terms in \eqref{eq: sync-gen-w} and \eqref{eq: sync-gen-i}) play the same role as $i_x$ and $v_x$. The self-synchronizing properties of a multi-machine power system are attributed to the exchange of kinetic and electrical energy through electrical torque and EMF pair. In the following section, we will assign this very mechanism for the inverter dynamics \eqref{eq: inverter dynamics}.

\section{Grid-Forming SM Matching Control} %
\label{sec: matching control}

From \cite{gross2017steady}, we know that every converter modulation controller inducing a synchronous, balanced, and sinusoidal steady-state must necessarily include an {\em internal model} of an oscillator of the form $\dot m_{\alpha\beta}^{\star} = \omega^{*}  {\boldsymbol{J}} m_{\alpha\beta}^\star$.
Thus, the first step in our design is to assign a sinusoidal modulation scheme parameterized in polar coordinates as
\begin{subequations}\label{eq: matching control}
\begin{align}
m_{\alpha\beta}(\theta) &=\mu \begin{bmatrix}
-\sin{\theta}\\ \cos{\theta}
\end{bmatrix}
\label{eq: internal model}
\,,
\end{align}%
where $\theta \in\mathbb{S}^1$ and {$\mu\in[0,1]$} are the {modulation's signal magnitude and angle, as controls} to be specified.
In the next step, we design a {grid-forming} modulation controller by {\em matching} the converter dynamics \eqref{eq: inverter dynamics}, augmented with the internal model \eqref{eq: internal model}, to the SM dynamics \eqref{eq: sync-gen}. By visual inspection we observe that model matching is achieved by {dynamic} feedback
\begin{equation}%
\dot\theta  = \eta v_{dc} 
\end{equation}%
\end{subequations}
where the constant $\eta = \omega_0/v_{dc,\textit{ref}}>0$  encodes the ratio between the nominal AC frequency $\omega_0$ and the DC voltage reference $v_{dc,\textit{ref}}$. All subsequent developments will be based on the matching control \eqref{eq: matching control}.

\begin{rem}{\bf(Equivalent SM interpretation)}
By defining the {equivalent angular velocity} as $\omega =\eta v_{dc}$ and by picking the modulation amplitude as $\mu=-2\eta L_m i_f$, we can rewrite $i_x$ and $v_x$ as
\begin{align}
\label{eq: closed loop i_x and v_x}
	i_{x} &= -\eta L_mi_f \begin{bmatrix} -\sin{\theta}\\ \cos{\theta} \end{bmatrix}^\top i_{\alpha\beta}
	\\
	v_{x} &= -L_mi_f \begin{bmatrix} -\sin{\theta}\\ \cos{\theta} \end{bmatrix} \omega \,.
\end{align}
We identify the AC-side switch voltage $v_x$ with the equivalent EMF voltage and the DC-side current $i_x/\eta$ with the {equivalent electrical torque in the machine.} %
Finally, we rewrite the closed loop \eqref{eq: inverter dynamics}, \eqref{eq: matching control} as the equivalent~SM
\begin{subequations}
\begin{align}
\dot \theta  &=\omega  
\\
\frac{C_{dc}}{\eta^2}\dot\omega  &= -\frac{G_{dc}}{\eta^2} \omega  + \frac{i_{dc}}{\eta} - \frac{1}{2\eta} m_{\alpha\beta}(\theta)^\top i_{\alpha\beta}
\\
L  \dot i_{\alpha\beta}  &= -R i_{\alpha\beta} - v_{\alpha\beta} + \frac{1}{2\eta} m_{\alpha\beta}(\theta)  \omega
\\
C \dot v_{\alpha\beta} &= -G v_{\alpha\beta}+ i_{\alpha\beta}-i_{l} 
\,,
\end{align}%
\label{eq: equivalent machine model}%
\end{subequations}%
where we identify $C_{dc}/\eta^2$, $G_{dc}/\eta^2$, and $i_{dc}/\eta$ with the equivalent mechanical inertia,
damping and mechanical driving torque, {respectively}. 
\oprocend
\end{rem}

\subsection{Closed-Loop Incremental Passivity}

In this section, we show how the matching controller \eqref{eq: matching control} can achieve desirable stability and passivity properties, while formulating them with respect to an induced operating point. 
Consider the closed-loop inverter dynamics \eqref{eq: inverter dynamics}, \eqref{eq: matching control}. By applying the $dq$-coordinate transformation \eqref{$dq$-transf} with angle $\theta $, we arrive at the following subsystem, which is independent of the angle state variable
\begin{subequations}%
\begin{align}%
C_{dc} \dot v_{dc} &=-G_{dc} v_{dc} + i_{dc} - \frac{\mu}{2} \mathbf{e}_2^\top i_{dq}
\\
L{\dot{i}_{dq}} &=-Ri_{dq} - {v_{dc}\eta L\boldsymbol{J}} i_{dq} + \frac{\mu}{2} \mathbf{e}_2 v_{dc} - v_{dq}
\\
C\dot v_{dq} &=-Gv_{dq} - {v_{dc}\eta C\boldsymbol{J}}v_{dq} - i_{l,dq} + i_{dq} \,,
\end{align}%
\label{eq: inv-dq}%
\end{subequations}%
where $i_{dq}\!=\!\boldsymbol{R}_\theta^\top i_{\alpha\beta}$, $v_{dq}\!=\!\boldsymbol{R}_\theta^\top v_{\alpha\beta}$ and $\frac{\mu}{2}\mathbf{e}_2\!=\!\boldsymbol{R}_\theta^\top m_{\alpha\beta}(\theta)$.

The following result characterizes the {strict} incremental passivity of the {$dq$-frame} inverter system \eqref{eq: inv-dq}, with respect to a steady-state {solution}, {as per Definition $1$ in \cite{trip2017passivity}}.

{\begin{thm}{\bf(Strict passivity in $dq$-frame)}
	Consider the model-matched system \eqref{eq: inv-dq} {and} assume that, for a given constant input $u^* = (i_{{dc}}^*, i_{l,dq}^*)$, there exists an equilibrium $x^{*} = (v_{{dc}}^{*},i_{dq}^{*},v_{dq}^{*})$ that satisfies
	\begin{equation} 
	\frac{C^2\norm{v_{dq}^*}^2}{4G} + \frac{L^2\norm{i_{dq}^*}^2}{4R}< \frac{G_{{dc}}}{\eta^2}
	\,.
	\label{eq: neg-def-cond}
	\end{equation}
	Then, system \eqref{eq: inv-dq} with input $ u = ( i_{{dc}}, - i_{l, dq})$ and output $ y = (v_{{dc}}, v_{dq})$ is strictly passive relative to the pair $(x^*, u^*)$.
\label{thm: closed-loop}  
\end{thm}}

\begin{pf}	
Our proof is inspired by \cite{YS-PT:14}. {Starting from the assumptions of the theorem, we define the} error coordinates $\tilde v_{dc}\!=\!v_{dc}\!-\!v_{dc}^*$, $\tilde i_{dq}\!=\!i_{dq}\!-\!i_{dq}^*$, $\tilde v_{dq}\!=\!v_{dq}\!-\!v_{dq}^*$, $\tilde{i}_{l,dq}\!=\!i_{l,dq}\!-\!i_{l,dq}^*$, $\tilde{i}_{dc}\!=\!i_{dc}\!-\!i_{dc}^{*}$, as well as $\omega^*\!=\!\eta v_{dc}^*$, {such that} the associated transient dynamics are expressed as
\begin{align}%
C_{dc}\dot{\tilde v}_{dc} =& -{G_{dc}}\tilde v_{dc} + \tilde i_{dc}-\frac{\mu}{2}\mathbf{e}_2^\top\tilde i_{dq}
\nonumber\\
L \dot {\tilde  i}_{dq} =&-{\left(R\boldsymbol{I}\!+\!{v_{dc}^*\eta L \boldsymbol{J} \!+\! \tilde{v}_{dc}\eta  L \boldsymbol{J}}\right)} \tilde  i_{dq} +\frac{\mu}{2} \mathbf{e}_2 \tilde v_{dc} 
\label{eq: sys-err-dq}
\\
&- \tilde v_{dc} \eta{L\boldsymbol{J}} i_{dq}^{*} - \tilde{v}_{dq}
\nonumber\\
C \dot {\tilde v}_{dq} =&-{\left(G\boldsymbol{I}\!+\!{v_{dc}^*\eta C\boldsymbol{J} \!+\! \tilde{v}_{dc}\eta  C \boldsymbol{J}}\right)} \tilde v_{dq}  - \tilde i_{l,dq}
\nonumber
\\
& - \tilde v_{dc} \eta C{\boldsymbol{J}}v_{dq}^{*}  + \tilde i_{dq}  \,.\nonumber
\end{align}
By considering the physical storage of the circuit elements, we define the incremental positive definite and differentiable storage function $\mathcal{V}_1:\real^{5}\to\real_{>0}$ as
\begin{equation}
\mathcal{V}_1=\frac{1}{2} C_{dc}\tilde v^2_{dc}+\frac{1}{2}L\tilde i_{dq}^{\top}\tilde i_{dq}+\frac{1}{2} C\tilde{v}_{dq}^{\top}\tilde{v}_{dq}
\label{eq: Lyapunov-func1}
\,.
\end{equation}
Due to the skew symmetry of ${\boldsymbol{J}}$, the derivative of $\mathcal{V}_1$ along the trajectories of the error system \eqref{eq: sys-err-dq} reads as
\begin{equation*}
\mathcal{\dot {V}}_1
= -
\begin{bmatrix} \tilde v_{dc} & \tilde i_{dq}^\top &\tilde v_{dq}^\top
\end{bmatrix} \mathcal{Q} \begin{bmatrix} \tilde v_{dc} & \tilde i_{dq}^\top &\tilde v_{dq}^\top
\end{bmatrix}^\top
- \tilde{v}_{dq}^{\top}\tilde{i}_{l,dq}  + \tilde{i}_{dc}\tilde v_{dc}\,,
\end{equation*}
where the symmetric matrix $\mathcal{Q} \in\real^{5\times 5}$ is given by
\begin{equation}
\mathcal{Q}  = \begin{bmatrix}
G_{dc} & \frac{1}{2}( \eta L\boldsymbol{J} i_{dq}^*)^{\top} & \frac{1}{2}(\eta C\boldsymbol{J} v_{dq}^*)^{\top}
\\
\frac{1}{2}(\eta L \boldsymbol{J} i_{dq}^*) & R\boldsymbol{I}& O
\\
\frac{1}{2}(\eta C\boldsymbol{J} v_{dq}^*)  & O & G\boldsymbol{I}
\end{bmatrix} 
\label{eq: P-matrix}
\end{equation}
By evaluating all leading principal minors of $\mathcal{Q}$ we see that under condition \eqref{eq: neg-def-cond}, $\mathcal{Q}$  {is} positive definite. Hence, system \eqref{eq: sys-err-dq} is strictly passive with input $(\tilde i_{dc}, -\tilde i_{l,dq})$ and output $(\tilde v_{dc}, \tilde v_{dq})$.  
\oprocend
\end{pf}

The importance of this result is that, when the load current $i_{l,dq}$ and the source current $i_{dc}$ are constant, the origin of \eqref{eq: sys-err-dq} is rendered asymptotically stable via Lyapunov's direct method.
Since ${\mathcal V}_{1}$ is radially unbounded, we obtain global asymptotic stability as well as the absence of any other type of equilibrium. 
{We shall further pursue this analysis after closing the passive ports of the inverter via a suitable DC actuation and an AC load current.}

\subsection{Closed-Loop Incremental Stability}

The {strict passivity} condition \eqref{eq: neg-def-cond} requires sufficiently large damping in the AC and DC components of the converter. However, parasitic resistances $R$ and $G_{dc}$ {can be arbitrarily} small in practice. To alleviate this shortage of stabilizing dissipation, we {implement} a DC-side {actuation} akin to governor speed droop control for generators to enforce condition \eqref{eq: neg-def-cond}. We propose the current source $i_{dc}$ to implement the proportional (P) controller %
\begin{equation}
i_{dc}=i_{dc,\textit{ref}} -K_{p} (v_{dc} - v_{dc,\textit{ref}})
\,,
\label{eq: idc-p-ctrl}
\end{equation}
with gain $K_p>0$ and set-points $i_{dc,\textit{ref}}>0$ and $v_{dc,\textit{ref}}>0$ for the DC-side current injection and the DC-side voltage, respectively.

{We are now ready to introduce the load model which we find best representative for the grid-forming application.} 
Assume that the load consists of a constant shunt impedance $\mathcal{Y}_l=G_l{\boldsymbol{I}}+ B_l\boldsymbol{J} $ accounting for passive devices (e.g., RLC circuits) connected to the converter. In parallel with this impedance, consider a sinusoidal current source with state $s_{l}$ having, for all time, the same frequency $\omega$ as the converter and otherwise constant amplitude. The latter can model a weak grid without grid-forming units, i.e., without any generator or inverter that regulates frequency and voltage, but possibly containing grid-following units equipped with PLLs which (instantaneously) synchronize to the frequency $\omega$.

\begin{assumption} The load current $i_l$ is given by the following system driven by the input $(\omega, v_{\alpha\beta})$,
		\begin{equation}
		\label{eq: load-model}
		\begin{split}
		\dot{s}_{l} &= \omega \boldsymbol{J}s_{l}  
		\\
		i_l &= (G_{l}\boldsymbol{I}+ B_l \boldsymbol{J}) v_{\alpha\beta} + s_{l} \,,
		\end{split}
		\end{equation}
		where $G_{l},B_l>0$ are constant parameters, and $s_{l}\in\real^2$ is the state of an internal oscillator. 
		\label{ass: iL-model}
	\end{assumption}
Notice that the internal state $s_{l}$ of this load model, when represented in the converter-side $dq$-coordinates, becomes a {constant}. 
All devices in the network can now be studied with respect to a single $dq$-frame angle, namely that of the virtual oscillator.
In this scenario we arrive at the following corollary. 

\begin{cor}{\bf(Closed-loop stability with DC-side P-control)} %
 Consider the inverter in {system} \eqref{eq: inv-dq} together with P-controller \eqref{eq: idc-p-ctrl} on the DC-side and load \eqref{eq: load-model} on the AC-side. Assume there exists a steady-state $x^{*} = (v_{dc}^{*},i_{dq}^{*},v_{dq}^{*})$ satisfying
\begin{equation}
\frac{C^2\norm{v_{dq}^{*}}^2}{4(G+G_{l})} + \frac{L^2\norm{i_{dq}^{*}}^2}{4R} < \frac{G_{dc}+K_p}{\eta^2}
\label{eq: neg-def-cond-PID}
\,.
\end{equation}
{Then,} for stationary {loads $s_{l, dq} = s_{l dq}^{*}$}, the steady state $x^{*}$ is unique and globally asymptotically stable. 
\label{cor: closed-loop plus DC P-control}
\end{cor}

Observe that condition \eqref{eq: neg-def-cond-PID} can be met by suitable choice of gain $K_p$ and that the condition is worst at no load, i.e. when $G_l = 0$. Furthermore, at this point, we cannot necessarily guarantee exact regulation of $v_{dc}$ to a particular $v_{dc,\textit{ref}}$ without having access to the load measurement. This discussion will be {addressed} later, in Section~\ref{sec: outer loops}.

Finally, the incremental passivity {property} highlighted in Theorem \ref{thm: closed-loop} is regarded as a key requirement for stability under interconnection, see \citep{FZO+13,YS-PT:14}, however this requires a single coordinate frame analysis for the networked scenario. Since in our work we use a $dq$-coordinate frame attached to a particular converter angle, the analysis does not pertain to a setup containing multiple (grid-forming) inverters. Nevertheless, this property is preserved in all our subsequent developments.

\subsection{Droop properties of matching control} %
\label{sec: properties}

An important aspect of {\em plug-and-play} operation in power systems  
is {steady-state} power sharing amongst multiple inverters by means of a droop characteristic. This is typically achieved via a trade-off between power injection and voltage amplitude or frequency \citep{FD-JWSP-FB:14a}. We {now} investigate these steady-state properties {which arise naturally in} the closed-loop system \eqref{eq: inverter dynamics}, \eqref{eq: matching control}, \eqref{eq: idc-p-ctrl}.

{Let ${r_x} =\frac{1}{2} \mu v_{dc}^*$ and ${\omega_x}=\eta v_{dc}^*$ denote the switching node voltage amplitude and frequency at steady state. Let ${P_x} = v_{x}^{\star\top} i_{\alpha\beta}^\star$ and ${Q_x} = v_{x}^{\star\top}\boldsymbol{J}^\top i_{\alpha\beta}^\star$ denote the active and reactive powers flowing from the switching node, at steady state, as per the convention in \citep{HA-YK-KF-AN:83}, and assume that they are constant.}

Two converters indexed by $i,j$ are said to achieve {\em proportional power sharing at ratio $\rho>0$} if
 $P_{x,i}/P_{x,j} = \rho$. %
{Furthermore}, the linear sensitivity factors relating steady-state active power injection ${P_x}$ to voltage amplitude ${r_x}$ and frequency ${\omega_x}$, are {defined} here as the {\em droop coefficients} ${d_r} = {\partial P_x}/{\partial {r_x}}$ and ${d_\omega} = {\partial P_x}/{\partial {\omega_x}}$. {Their relationship is given in the proposition that follows.}

\begin{prop}{\bf(Droop slopes)} 
Consider {system} \eqref{eq: inverter dynamics}, together with matching controller \eqref{eq: matching control} and DC-side controller \eqref{eq: idc-p-ctrl}. {Denote} $ {\omega_x}=\omega ^*$, and define {the constant} $i_{0}= i_{dc,\textit{ref}}+ K_p v_{dc,\textit{ref}}$. The following statements hold at equilibrium: 
\begin{enumerate}
\item Nose curves: the switching voltage amplitude {${r_x}$ has} the following expression as a function of {$i_0$ and $P_x$}
\begin{align*}
{r_x} &=\frac{\mu}{4(G_{dc}+K_p)}{\left(i_{0} {\pm}\sqrt{{i_{0}}^{2}-4(G_{dc}+K_p)P_x}\right)},
\end{align*}%
with a similar expression for the virtual frequency ${\omega_x} = \tfrac{2\eta}{\mu}{r_x}$.
Moreover, the reactive power ${Q_x}$ {and the quantities ${r_x},{\omega_x}$ are not related}.

\item  Droop behavior: {around the} operating point ${\omega_x}$, the expression for the frequency droop coefficient is given by
\begin{equation}
\label{eq: droop slopes}
	{d_\omega}=-\frac{2(G_{dc}+K_p)}{\eta^2}{\omega_x}+\frac{i_{0}}{\eta} \,,
\end{equation}
with {an analogous} expression for the switching node voltage amplitude droop $d_r$, since ${r_x} = \tfrac{\mu}{2\eta}{\omega_x}$.
\item Power sharing:  Consider a pair of converters $i$ and $j,\,\{i,j\}\in \mathbb{N}$ with $G_{dc}=0$, {identical DC-side} voltage references $v_{dc,\textit{ref}}>0$ and control gain $\eta>0$. The converters achieve proportional power sharing at ratio $\rho = P_{x,i}/P_{x,j}$ if 
\begin{align} 
K_{p,i} ={\rho} K_{p,j} \; , \; i_{dc,\textit{ref},i} = {\rho} i_{dc,\textit{ref},j} \,,
\label{eq: pow-sharing}
\end{align}
or equivalently if $d_{\omega,i}=\rho {d_{\omega,j}}$ and ${P_{dc,i}}=  \rho {P_{dc,j}}$ {with $P_{dc}=v_{dc,\textit{ref}}\cdot i_{dc,\textit{ref}}$}\,.
\end{enumerate}
\label{prop: ss-charcterization}
\end{prop}

\begin{pf} 
   To prove statement $(1)$ consider the DC{-side} dynamics \eqref{eq: DC voltage} at steady state
\begin{equation}
	0=-(G_{dc}+K_p) v_{dc}^*+i_{0} - i_x^* \,.
   \label{eq: steady state vdc}
\end{equation}
We multiply \eqref{eq: steady state vdc} by $v_{dc}^*$ to obtain quadratic expression relating $P_x=i_x^*\cdot v_{dc}^*$ and $v_{dc}^*$, at steady state
\begin{equation}
v_{dc,\pm}^*=\frac{i_{0} {\pm}\sqrt{i_{0}^{2} -4(G_{dc}+K_p)P_x}}{2 (G_{dc}+K_p)}
\,.
\label{eq: v_dc at ss}
\end{equation}
The claimed nose curves follow directly. Consider now
 \begin{align}
 \label{eq: power equations}
 {P_x} &= -\frac{(G_{dc}+K_p)}{\eta^{2}} \omega^{2}_x+ \frac{i_{0}}{\eta} {\omega_x} \,.
\end{align}
By linearizing the above equations around the steady-state operating point ${\omega_x}$, we find the droop slope in \eqref{eq: droop slopes}.
Finally, the proportional power sharing ratio $\rho>0$ between two converters $i$ and $j$ is given by setting $G_{dc} = 0$ in \eqref{eq: power equations}
\begin{equation}
\rho = \frac{P_{x,i}}{P_{x,j}}=\frac{\frac{K_{p,i}}{\eta^{2}} {\omega_x}- \frac{i_{0,i}}{\eta}}{\frac{K_{p,j}}{\eta^{2}} {\omega_x}- \frac{i_{0,j}}{\eta}}\,.
\label{eq: ac power}
\end{equation}
The latter equality is satisfied if \eqref{eq: pow-sharing} holds. 
\oprocend
\end{pf} 
{The following remarks can now be drawn: Statement \textit{(1)} gives two solutions\footnote{for the rest of this paper our system is not subject to constant power loads and has an unique admissible equilibrium} 
 for the voltage amplitude $r_{x}$. Among these two, the so-called {\em high-voltage solution} (with the plus sign) is the practically relevant operating point as depicted in Figure \ref{Fig: droopslopes}. From statement \textit{(1)}, we can also deduce that the maximal active power which can be delivered at the switching node, $P_{max} = {i_{0}^{2}}/ \left({4(G_{dc}+K_p)}\right)$, is marked by the right tip of the nose curve. No stationary solutions exist beyond this bifurcation point.
\begin{figure}[h!]
	\centering
	\includegraphics[width=\columnwidth]{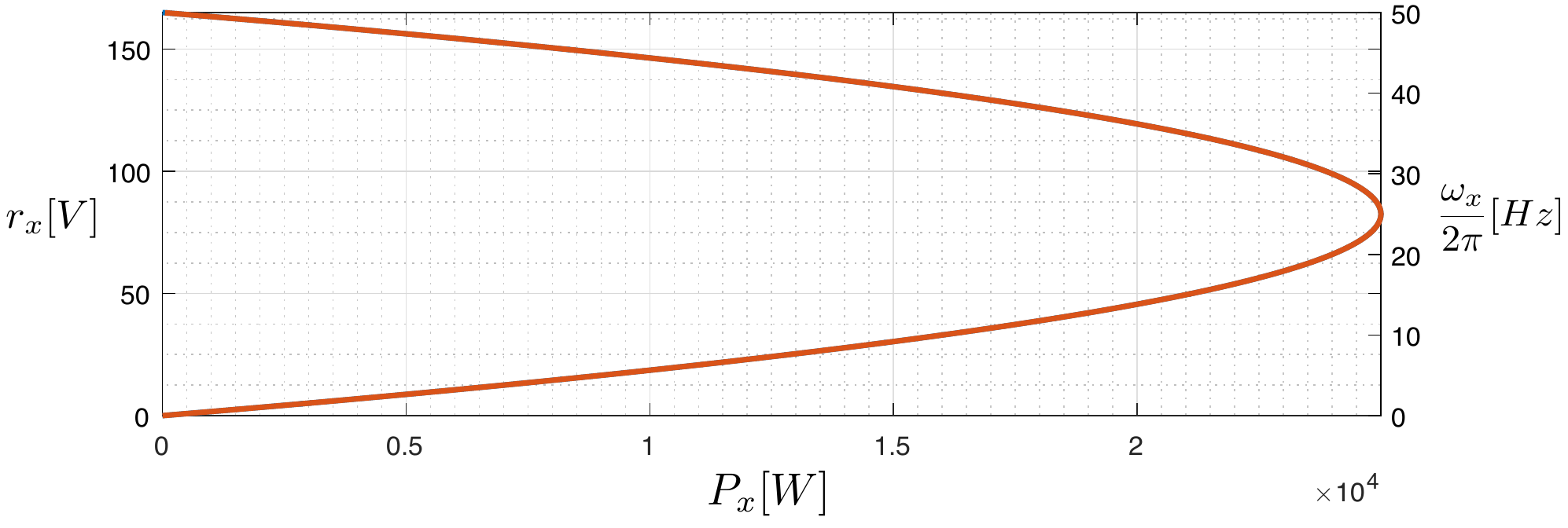}
	\caption{Steady-state {profiles} $({r_x},{P_x})$ and $({\omega_x},{P_x})$ for the set of converter parameters described in Section \ref{sec: case study}.}
	\label{Fig: droopslopes}
\end{figure}
To the best of our knowledge, a typical inverter design will have, by design, its operating region away from the tip of the nose curve, where the linear sensitivity factors are a good approximation and parametric bifurcations are of no practical concern.} 
Regarding statement {\it(3)}, the power sharing conditions \eqref{eq: pow-sharing} are perfectly analogous to the ones in conventional droop control \citep{FD-JWSP-FB:14a}: the droop slopes and the power set-points must {be related by} the same ratio $\rho$.  
Finally, we remark that similar expressions as in \eqref{eq: pow-sharing} can be obtained for a non-zero DC-side damping and heterogeneous converter parameters.%

\subsection{Relation to Other Converter Control Strategies} 

Our matching control  can be understood from the viewpoint of PBC by writing the inverter \eqref{eq: inverter dynamics} as the Port-Hamiltonian system \citep{AJvdS:96}
\[
\dot z =\left(\mathcal{J}(m)-\mathcal{D}\right)\nabla H(z) +  \mathcal{G}u \,,
\]
where $z = (C_{dc}v_{dc}, L i_{\alpha\beta}, Cv_{\alpha\beta})$ is the state, $m$ is the modulation, $u = (i_{dc}, -i_{l})$ is an exogenous input, $H(z) = \frac{1}{2} C_{dc}v_{dc}^{2}$ $+ \frac{1}{2}i_{\alpha\beta}^\top Li_{\alpha\beta}$ $+ \frac{1}{2} v_{\alpha\beta}^\top Cv_{\alpha\beta}$ is the physical energy, as in \eqref{eq: Lyapunov-func1}, $\mathcal{J}(m)$ is {a} skew-symmetric interconnection matrix {depending on the modulation signal $m$}, $\mathcal{D}$ and $\mathcal{G}$ are positive definite damping and input matrices. 
The Port-Hamiltonian structure is preserved upon augmenting the inverter with the internal model \eqref{eq: internal model}. On this ground, we can link our approach to that of PBC and IDA-based matching control \citep{ortega2004interconnection}. In particular, matching controller \eqref{eq: matching control} together with P-controller \eqref{eq: idc-p-ctrl} can be understood as IDA reshaping the ${ \mathcal J}$ and ${ \mathcal D}$ matrices.

Our control strategy {can also be associated with} oscillator-based controller methods.	By defining $m\in\real^2$ as the controller state, we can rewrite \eqref{eq: matching control} as %
\begin{equation*}
\dot m = \omega \boldsymbol{J} m  \,,
\label{eq: controller dynamics}
\end{equation*}
i.e, the matching control \eqref{eq: matching control} is an oscillator with constant amplitude $\norm{m(0)}=\mu$ and state-dependent frequency $\omega =\eta v_{dc}$ as feedback for the converter dynamics \eqref{eq: inverter dynamics}. This control strategy  resembles the classic {\em proportional resonant control} \citep{RT-FB-LM-LPC:06} with the difference that the frequency of the oscillator \eqref{eq: controller dynamics} adapts to the DC voltage which again reflects the grid state. 
Another related control strategy is {\em virtual oscillator control}  encoding the inverter terminal dynamics as a nonlinear limit cycle oscillator adapting to the grid state \citep{BBJ-SVD-AOH-PTK:14a,sinha2017uncovering}.
\section{Voltage and frequency regulation}
\label{sec: outer loops}


Starting from the model-matching controller \eqref{eq: matching control}, we now look to design outer control loops for the current source $i_{dc}$ as well as the modulation amplitude $\mu$ with the aim of tracking a given reference, initially for the DC capacitor voltage and then also for the AC capacitor voltage amplitude.

\subsection{Exact Frequency Regulation via Integral Control}
\label{subsec: Frequency Regulation via Virtual Governor Control}

{In some scenarios, e.g., in islanded microgrids, it is desirable that inverters also contribute to frequency regulation (usually called secondary control) rather than mere droop control.}
Inspired by frequency regulation of SMs via governor control, i.e., controlling the torque in \eqref{eq: sync-gen} as a function of the frequency, we propose a frequency regulation strategy by pairing the passive input and output, $\tilde i_{dc} = i_{dc}-i_{dc,\textit{ref}}$ and $\tilde v_{dc} = v_{dc}-v_{dc,\textit{ref}}$ respectively, in the inverter model \eqref{eq: inv-dq}, in the aim of tracking a reference frequency $\omega_0=\eta v_{dc,\textit{ref}}$. We {propose} the PID controller%
\begin{equation} 
i_{dc}
=i_{dc,\textit{ref}}-K_{p}\tilde v_{dc}-K_{i} {\int_{0}^{t} \tilde {v}_{dc}(\tau)\,d\tau} -K_{d}\dot {\tilde v}_{dc}
\,,
\label{eq: idc-ctrl}
\end{equation}
where $i_{dc,\textit{ref}}>0$ {is a user-defined parameter,} and $K_{p}$, $K_{i}$, $K_d$ are positive control gains.
{Before proceeding to the stability result we put together system \eqref{eq: inv-dq}, controller \eqref{eq: idc-ctrl}, and load model \eqref{eq: load-model}, and express the closed loop in error coordinates formulated relative to an induced equilibrium. To account for the newly introduced integral term, we define the state variable $\xi\in\real$ and denote its steady-state value by $\xi^{*}$ such that $\tilde\xi = \xi - \xi^{*}$.}
\begin{align}
\dot {\tilde \xi} =&~{\tilde v_{dc}}
\nonumber\\
(C_{dc}+{K_d}) \dot{\tilde v}_{dc}=&-(G_{dc}+K_{p}) \tilde v_{dc} - K_{i} \tilde \xi
-\frac{\mu}{2}\mathbf{e}_2^\top\!\tilde i_{dq}
\nonumber\\
L \dot {\tilde  i}_{dq} =&-{ \left(\mathcal{Z}_{} + \tilde v_{dc} \eta L \boldsymbol{J}\right)} \tilde i_{dq} + \frac{\mu}{2} \mathbf{e}_2{\tilde v}_{dc}
\nonumber\\ 
\label{eq: dq-idc-ctrl}%
& - \tilde v_{dc} \eta L\boldsymbol{J} i_{dq}^{*} - {\tilde v}_{dq} 
\\
C \dot {\tilde v}_{dq} =& -{\left(\mathcal{Y}_{} +\mathcal{Y}_l + \tilde v_{dc} \eta C \boldsymbol{J}\right)} \tilde v_{dq} \nonumber \\
&- \tilde v_{dc} \eta C\boldsymbol{J}v_{dq}^{*} + \tilde i_{dq} \nonumber
\end{align}%
%
where $\mathcal{Z}_{} = (R\boldsymbol{I}+\omega_0L\boldsymbol{J})$ and $\mathcal{Y}_{} = (G\boldsymbol{I}+\omega_0C\boldsymbol{J})$ 
are the impedance of the AC-side inductor and the admittance of the AC-side capacitor, respectively. Since PID control of the DC-voltage \eqref{eq: idc-ctrl} is common practice in DC/AC converters, {we will see that} pairing it with the matching control \eqref{eq: matching control} yields {exact} AC frequency regulation. 

The following result addresses existence, uniqueness, and stability of a desired steady state of the closed-loop system \eqref{eq: dq-idc-ctrl} satisfying $v_{dc}^*\!=\!v_{dc,\textit{ref}}$ and $\omega^*\!=\!\omega_0$. Typically, in such systems, $\omega_0$ can be seen as the grid nominal frequency, while $v_{dc,\textit{ref}}$ the {reference} voltage of the converter's DC-link capacitor. By appropriately choosing the gain $\eta=\omega_0/v_{dc,\textit{ref}}$, we are able to achieve both specifications.

\begin{thm}[Exact frequency regulation] %
Consider the closed-loop system \eqref{eq: dq-idc-ctrl} {and a given set-point $\omega_0>0$}. The following {two} statements hold:
\begin{enumerate}
\item There exists a unique steady state at the origin with $\omega^* = \omega_0$ as system frequency.
  \item {Assuming} condition \eqref{eq: neg-def-cond-PID} is satisfied, the zero-equilibrium of \eqref{eq: dq-idc-ctrl} is globally asymptotically stable. 
\end{enumerate}
\label{thm: governor-ctrl}
\end{thm}

\begin{pf}
A {steady state} of the closed loop \eqref{eq: dq-idc-ctrl} is characterized by $\tilde v_{dc}=0$ and a linear set of equations $A \left[\begin{smallmatrix}
\tilde \xi & \tilde i_{dq}^\top & \tilde  v_{dq}^\top \end{smallmatrix}\right]^{\top}=0$,
where $A\in\real^{5\times5}$ is given by
\begin{equation}
\label{eq: A-matrix}
A = \begin{bmatrix}	-K_i & -\tfrac{\mu}{2}\mathbf{e}_2^\top & O \\
O & -\mathcal{Z}_{} & \boldsymbol{I} \\
O & -\boldsymbol{I} & -(\mathcal{Y}_{}+\mathcal{Y}_l) \end{bmatrix}  \in \real^{5 \times 5} \,.
\end{equation}
It follows that
$$\det(A) = -K_i \norm{\mathcal{Z}_{} (\mathcal{Y}_{}+\mathcal{Y}_l)+\boldsymbol{I}}^2\,,$$
such that positivity of the converter and load parameters assures {invertibility} of $A$, and hence
$\left[\begin{smallmatrix}
\tilde \xi & \tilde i_{dq}^\top & \tilde  v_{dq}^\top
\end{smallmatrix}\right]^{\top}=0$. 
 Thus, there is a unique zero steady-state for this error subsystem. The stability proof of this steady state is analogous to the proof of Corollary \ref{cor: closed-loop plus DC P-control}
after replacing the original storage function $\mathcal{V}_1$ with 
$
 \mathcal{V}_{2}=   \mathcal{V}_{1}+  \frac{1}{2}K_{i} \tilde \xi^{2}+  \frac{1}{2} K_d \tilde v_{dc}^2,\,
$
to account for $\tilde \xi$ {and the {gain} $K_{d}$.}
 With these modifications the derivative of the storage function $\mathcal{ V}_2$ becomes %
 \begin{equation*}
\mathcal{\dot V}_{2}= -
\begin{bmatrix} \tilde v_{dc} & \tilde i_{dq}^\top & \tilde v_{dq}^\top
\end{bmatrix} \mathcal{Q}  \begin{bmatrix} \tilde v_{dc} & \tilde i_{dq}^\top &\tilde v_{dq}^\top
\end{bmatrix}^\top \leq 0\,,
\end{equation*}
where $\mathcal{Q}$ is as in \eqref{eq: P-matrix} with $G_{dc}$ {and $G$} replaced by $G_{dc}+K_{p}$ {and $G+G_{l}$, respectively.} Finally a LaSalle{-type} argument accounting for the state $\tilde \xi$ together with radial unboundedness of $\mathcal{V}_{2}$ 
guarantees global asymptotic stability.
\oprocend
\end{pf}
{Notice that} the P-control on the DC voltage enhances the overall system stability, as discussed before. Furthermore, by comparing systems 
\eqref{eq: inv-dq} 
and \eqref{eq: dq-idc-ctrl}, observe that the effect of the {PID} gains is to provide additional {inertia and damping to the DC circuit}. 
Lastly, from a conventional power system perspective, it is instructive to write the frequency error dynamics, whereby $\tilde\omega = \eta \tilde{v}_{dc}$  
\begin{equation*}
\frac{(C_{dc}+ K_d)}{\eta^2} \dot{\tilde\omega} = -\frac{(G_{dc}+ K_{p})}{\eta^2}\tilde \omega -\frac{K_{i}}{{\eta}} \int_0^t \tilde \omega(\tau) d\tau  -\frac{1}{\eta}i_x
\,,
\end{equation*}
which for $K_{i}=0$ resemble the standard swing equations with synthetic droop and inertia induced by $K_{p}$ and $K_{d}$.

We conclude that for secondary frequency regulation {-- independently} of the particular {modulation} strategy -- a sufficiently large equivalent DC energy storage is required to cope with {a given power} imbalance. If the task of frequency regulation is to be shouldered by multiple inverters, then the decentralized integral control in \eqref{eq: idc-ctrl} can be easily adapted to broadcast AGC-like or consensus-based distributed integral control schemes \citep{de2018bregman,FD-JWSP-FB:14a,FD-SG:16}, which assure robust power sharing.

\subsection{Amplitude  Regulation by {Disturbance Feedback}}
This section investigates a series of controllers designed to regulate the AC-side voltage amplitude $\norm{v_{\alpha\beta}}$ to a desired set-point $r_{\textit{ref}}>0$. Throughout this section, we assume that the load in Assumption \ref{ass: iL-model} has zero shunt impedance, namely that $\mathcal{Y}_l=0$, i.e., the load is purely of constant- (in $dq$-frame) current nature, so that $i_{l,dq}=s_{l,dq}$. 
This modeling choice is not merely done for simplicity of exposition (the load impedance can always be absorbed in the filter conductance $G$) but mainly due to the fact that all amplitude controllers (most importantly, droop control) explicitly or implicitly rely on a measurement of the load current which is considered to be an exogenous signal.

\subsubsection{Feasibility and feedforward control} 

We now consider as actuation input, the modulation amplitude $\mu$, analogously to standard practice in SM excitation current. Formally our control objective is to achieve $\norm{v_{dq}}=r_{\textit{ref}}$, at steady state. Let us first characterize the feasibility of this task in terms of the system parameters. 

\begin{thm}[Existence of {load-induced} equilibria]
	Consider the closed-loop inverter model \eqref{eq: dq-idc-ctrl} with $\mathcal{Y}_l=0$.
	For given set-points $r_{{\textit{ref}}}>0,\, v_{{dc, \textit{ref}}}>0$ and constant load current $s_{l,dq} \in \real^2$, define  the quantity
	\begin{equation}
	\psi = r_{{\textit{ref}}}^2 \norm{\mathcal{Z}_{}\mathcal{Y}_{}+\boldsymbol{I}}^2 - \norm{\mathcal{Z}_{} s_{l,dq}}^2
	\label{eq: product}
	\,,
	\end{equation}
	Then, the following statements are equivalent:
	\begin{enumerate}
		\item There exists a unique steady state $(\xi^*,v_{{dc}}^*, i_{dq}^*, v_{dq}^*)$ that satisfies $\norm{v_{dq}^*}=r_{{\textit{ref}}}$ and $\mu>0$; and
		\item $\psi > 0$.
	\end{enumerate}
	\label{thm: eq-existence}
\end{thm}

\begin{pf}	
	We formulate the equilibria of {system} \eqref{eq: dq-idc-ctrl}, together with the requirement that {$\norm{v_{dq}^*} = r_{{\textit{ref}}}$}, as
	\begin{subequations}
		\begin{align}
		0 &= v_{{dc}}^*-v_{{dc},{\textit{ref}}}  
		\label{eq: ss-eq-1}%
		\\
		0 &=-(G_{{dc}}+K_p) v_{{dc}}^* -K_{i}\xi^* - \frac{\mu}{2} \mathbf{e}_2^{\top} i_{dq}^* 
		\label{eq: ss-eq-2}\\
		0&=-{(R\boldsymbol{I}+v_{dc}^*\eta L{ \boldsymbol{J}})} i_{dq}^* + \frac{\mu}{2} \mathbf{e}_2 v_{{dc}}^*-v_{dq}^*
		\label{eq: ss-eq-3}\\%
		0 &=-{({G\boldsymbol{I}+v_{dc}^*\eta C}{\boldsymbol{J}} )} v_{dq} ^*- s_{l, dq} + i_{dq}^*
		\label{eq: ss-eq-4}\\
		0&=v_{dq}^{*\top}v_{dq}^*-r_{{\textit{ref}}}^2,\,
		\label{eq: ss-eq-5}%
		\end{align}%
		\label{eq: des-ss-eq}%
	\end{subequations}%
	By subsequent elimination of variables, we can solve equations \eqref{eq: des-ss-eq} for $\mu$ in terms of input $s_{l,dq}$ and set-point $v_{dc,ref}$. We arrive at the quadratic equation 
	\begin{subequations}
		\begin{align*}
		0&=\mu^2-{b}\mu- \frac{4\psi}{v_{{dc},{\textit{ref}}}^2} \,,
		\end{align*}
	\end{subequations}
	{where ${b}=\frac{4}{v_{{dc},{\textit{ref}}}}\mathbf{e}_2^\top \mathcal{Z}_{}^\top s_{l,dq}$} is the sum of the two solutions $\mu_{\pm}$ of the quadratic equation.
	These  solutions
	\begin{align}%
	\mu_{\pm} = \frac{{b}}{2}\pm \sqrt{\left(\frac{b}{2}\right)^2+ \frac{4\psi}{v_{{dc},{\textit{ref}}}^2}}
	\label{eq:mu-root}
	\end{align}%
	are real-valued and have opposite signs $\mu_+ > 0$, $\mu_- < 0$ if and only if $\psi > 0$.
	In what follows, we restrict ourselves to the unique positive solution $\mu_{+}>0$. 
	Notice from \eqref{eq: ss-eq-1} that $v_{{dc}}^*=v_{{dc}, {\textit{ref}}}$.  After replacing $\mu_{+}$ into \eqref{eq: ss-eq-2}-\eqref{eq: ss-eq-4}, the remaining equations are linear $A  \left[\begin{smallmatrix}\tilde \xi & \tilde i_{dq}^\top & \tilde  v_{dq}^\top \end{smallmatrix}\right]^{\top}=0$ with $A$ nonsingular, as in the proof of Theorem \ref{thm: governor-ctrl}. These equations can be solved uniquely for $( \xi^*, i_{dq}^*, v_{dq}^*)$ which is consistent with \eqref{eq: ss-eq-5} by choice of $\mu_{+}>0.$
	\oprocend \end{pf} 

{The condition $\psi>0$ can be interpreted as an upper bound for the admissible constant load current $s_{l,{dq}}$ as a function of the given set-point $r_{\textit{ref}}$, since otherwise there would be no real-valued solution for $\mu$. %
Observe that the constraint $\mu_+ \leq 1$ can be enforced by adjusting the converter parameters and by further limiting the maximum allowable load.}

\paragraph*{Disturbance decoupling control}
Starting from the insights given by Theorem~\ref{thm: eq-existence}, we are able to construct a disturbance-feedback, asymptotic output tracking controller which relies on measurement of the load current $s_{l,dq}$ to produce the modulation input $\mu$ according to \eqref{eq:mu-root}. This approach can be regarded as a system inversion of the transfer path from $s_{l,dq}$ to the regulated voltage output $\|v_{dq}\|$, a standard procedure in disturbance decoupling. In the next subsection, we will discuss two extensions to this control strategy following PBC and droop control specifications.

\begin{cor}[Disturbance decoupling control]%
\label{cor: disturbance decoupling}%
Consider {system} \eqref{eq: dq-idc-ctrl} with $\mathcal{Y}_l=0$.  {Assume that the} load disturbance $i_{l,dq} = s_{l,dq}$ is a constant and measurable signal and that $\psi$ defined in \eqref{eq: product} is positive. Given a reference AC voltage $r_{\textit{ref}}>0$, assign the modulation amplitude
$$\mu = \mu_+(s_{l,dq}),\,\nonumber$$

with $\mu_+$ is as in \eqref{eq:mu-root}.
Further assume that {passivity condition} \eqref{eq: neg-def-cond} holds. Then the unique equilibrium is characterized by $v_{dc}^* = v_{dc,\textit{ref}}$ and $\norm{v_{dq}^*} = r_{\textit{ref}}$, and is globally asymptotically stable for the closed loop.
\end{cor}

\begin{pf}
For any {constant} $\mu>0$ the desired closed-loop equilibrium is described by \eqref{eq: des-ss-eq}. The existence of such equilibrium is guaranteed under the condition $\psi>0$
and for $\mu_{\pm}$ as in \eqref{eq:mu-root}. By assigning the positive solution in \eqref{eq:mu-root}, the amplitude $\mu = \mu_{+}(s_{l,dq})$ is constant for a {given} constant load $s_{l,dq}$. 
The stability claim now follows from {the same} reasoning as in the proof of Theorem \ref{thm: governor-ctrl}.
\oprocend
\end{pf}

\subsubsection{Compatibility with existing control techniques}

While very effective in achieving the prescribed steady-state specification, notice that the {disturbance feedback} control in Corollary~\ref{cor: disturbance decoupling} requires exact knowledge of the plant as well as load measurement. 
To assess the robustness of this framework, we investigate two additional extensions which could provide some insight into the practicality of the implementation. 

\paragraph*{PI-PBC:}
{Inspired by \cite{zonetti2014globally}, we now derive a PI-PBC feedback by identifying a the passive output corresponding to the new considered input.} In this regard, with pick input $\mu = \tilde\mu + \mu^*$, with $\mu^*=\mu_+$ from \eqref{eq:mu-root} and $\tilde\mu$ yet to be designed. We we rewrite \eqref{eq: dq-idc-ctrl} with $\mathcal{Y}_l=0$ as 
%
\begin{align} \label{eq: sys-non-constant-mu}
\dot{\tilde \xi} =& ~\tilde v_{dc} 
\nonumber\\
(C_{dc}\!+\!{K_d}) \dot{\tilde v}_{dc} =& -(G_{dc}+K_p)\tilde v_{dc} -K_i\tilde \xi \nonumber \\
&-\frac{\tilde\mu + \mu_+}{2}\mathbf{e}_2^{\top} \tilde i_{dq}- \frac{\tilde\mu}{2}\mathbf{e}_2^{\top} i_{dq}^{*}
\\
L \dot{\tilde i}_{dq} =&-{\left(\mathcal{Z}_{} + \tilde{v}_{dc}\eta L\boldsymbol{J}\right)} \tilde i_{dq}\!-\!\tilde{v}_{dc}\eta L{\boldsymbol{J}} i_{dq}^*- \tilde{v}_{dq}
\nonumber
\\
&+ \frac{\tilde \mu+\mu_+}{2}\mathbf{e}_2 \tilde v_{dc} + \frac{\tilde \mu}{2}\mathbf{e}_2  v_{dc}^{*}
\nonumber
\\
C \dot{\tilde{v}}_{dq} =&-{\left(\mathcal{Y}_{} + \tilde{v}_{dc}\eta C\boldsymbol{J} \right)} \tilde {v}_{dq} - \tilde{v}_{dc}\eta C{ \boldsymbol{J}} v_{dq}^* + \tilde i_{dq} . 
\nonumber
\end{align}
For the computiation of $\mu_+$ in \eqref{eq:mu-root} we also assume a constant measurable load current $i_{l,dq}$, such that
the prescribed equilibrium of \eqref{eq: sys-non-constant-mu} satisfies \eqref{eq: des-ss-eq}, or equivalently $\psi>0$. Observe that system \eqref{eq: sys-non-constant-mu} is {passive} with input $\tilde \mu$, output $y= \tilde i_q v_{dc}^{*}- i_q^* \tilde v_{dc}$ and storage function ${\mathcal{V}}_2$ from before, since%
$$\dot {{\mathcal V}}_2= -\begin{bmatrix} \tilde v_{dc}& \tilde i_{dq}^\top &\tilde v_{dq}^\top
\end{bmatrix} \mathcal{Q}  \begin{bmatrix} \tilde v_{dc}& \tilde i_{dq}^\top &\tilde v_{dq}^\top
\end{bmatrix}^\top + \tilde \mu^{\top} { y}\,.$$
This last observation motivates the PI-PBC feedback
\begin{subequations}
\label{eq: amp-PBCc-ctrl}
\begin{align}%
\dot {\tilde\nu}&={ y} 
\\
{\tilde\mu}&= -\kappa_p { y} - \kappa_i {\tilde\nu} \,, %
\end{align}%
\end{subequations}
where ${y}= \tilde i_q v_{dc}^{*}- i_q^* \tilde v_{dc}$ and $\kappa_p, \kappa_i>0$. {Finally, the resulting error feedback becomes 
$$y= i_q v_{dc,\textit{ref}} - i_q^*v_{dc}\,,$$
with $i_q^* = \mathbf{e}_2^\top(\mathcal{Z}_{}+\mathcal{Y}_{}^{-1})^{-1}(\tfrac{\mu}{2}\mathbf{e}_2v_{dc,\textit{ref}}+\mathcal{Y}_{}^{-1}s_{l,dq})$.
This is the same {type} of output to be regulated to zero as identified in \citep{zonetti2014globally} indicating a power imbalance across the inverter.} %
\begin{prop}{\bf(PI-PBC)} %
Consider {system} \eqref{eq: sys-non-constant-mu} {with the} PI-PBC {feedback} \eqref{eq: amp-PBCc-ctrl}.  Assume that {the} load disturbance $s_{l,dq}$ is a constant measurable signal and that $\psi$ defined in  \eqref{eq: product} is positive. Further assume that the passivity condition \eqref{eq: neg-def-cond} holds. Then the unique equilibrium is characterized by $v_{dc}^* = v_{dc,\textit{ref}}$ and $\norm{v_{dq}^*} = r_{\textit{ref}}$, and is globally asymptotically stable.
\end{prop}

\begin{pf}
Consider the radially unbounded Lyapunov function
$\mathcal{ V}_3 = \mathcal{{V}}_2 + \frac{\kappa_i}{2}  \tilde \nu^2$ and its derivative along trajectories of  \eqref{eq: sys-non-constant-mu}, \eqref{eq: amp-PBCc-ctrl}
\begin{equation*}
\mathcal{\dot {V}}_3 
=
- \begin{bmatrix} \tilde v_{dc}& \tilde i_{dq}^\top &\tilde v_{dq}^\top 
\end{bmatrix} \mathcal{Q} \begin{bmatrix} \tilde v_{dc}& \tilde i_{dq}^\top &\tilde v_{dq}^\top
\end{bmatrix}^\top
-\kappa_p {y}^2
\leq 0\,,
\end{equation*}
where $\mathcal{Q}$ is as in \eqref{eq: P-matrix} with $G_{dc}$ replaced by $G_{dc}+K_{p}$.
Assuming condition \eqref{eq: neg-def-cond-PID} is met, a LaSalle argument accounting for the evolution of $\tilde \xi$ and $\tilde \nu$ guarantees global asymptotic stability.
\oprocend
\end{pf}

{Observe that the PI-PBC strategy \eqref{eq: amp-PBCc-ctrl} requires that the load current $i_{l,dq} = s_{l,dq}$ is measurable, as it is used in the computation of the steady-state inducing terms $\mu_+$ and $i_{q}^{*}$ in the feedback law. In this way, the feedforward control \eqref{eq:mu-root}, as well as PI-PBC \eqref{eq: amp-PBCc-ctrl}, endow the closed-loop with the ability of rejecting the disturbance $i_{l,dq}$ as long as it is constant or, due to global asymptotic stability, provided it eventually settles to a constant which is state-independent.}

\paragraph*{{Voltage} droop control:} 
We have seen that, without an integral term, the matching control has the inherent droop properties of the SM. 
However, by introducing controller \eqref{eq: idc-ctrl} for exact frequency regulation, this droop effect has been removed in both the amplitude and the frequency of the AC-side voltage. 
In the remainder of this subsection, we propose a control strategy that implements, instead, a voltage-power droop behavior. {The droop dependency will be based on active power measurement, nevertheless the same reasoning applies for the case of reactive power.} We start by choosing
\begin{equation}
\mu= \mu_{\textit{ref}}+d_v (P_{l}-P_{\textit{ref}}),\,\,
\label{eq: amp-droop-ctrl}
\end{equation} 
where $\mu_{\textit{ref}}=\tfrac{2r_\textit{ref}}{v_{dc,\textit{ref}}}$. Here, $r_\textit{ref}$ and $P_{\textit{ref}}$ are set-points for the AC voltage amplitude and load power, respectively, $d_v>0$ is the droop coefficient, and $P_l = i_{l,dq}^{\top} v_{dq}$ denotes the entire load power measurement. {The droop factor $d_v$ represents a linear trade-off} between the modulation amplitude $\mu$ and the active power ${P_{l}}$ and induces a steady-state amplitude $\norm{v_{dq}^*}$ that is not necessarily equal to the {prescribed} reference $r_{\textit{ref}}$. 
The aim of the following result is to show that {this particular droop strategy is also} compatible with our framework. 
\begin{prop}[{Voltage} droop control]
\label{prop: amp-droop-ctrl}
Consider {system} \eqref{eq: dq-idc-ctrl} {with input $\mu$ given by \eqref{eq: amp-droop-ctrl} and $\mathcal{Y}_l=0$. 
Further assume that the closed loop \eqref{eq: dq-idc-ctrl}, \eqref{eq: amp-droop-ctrl} admits a steady state $(\tilde\xi, \tilde v_{dc},\tilde i_{dq},\tilde v_{dq}) = 0$.} Assuming that condition \eqref{eq: neg-def-cond} holds, then for a sufficiently small droop coefficient $d_v>0$, this steady state is globally asymptotically stable for the closed loop. 
\end{prop}

\begin{pf}
 We  rewrite the closed-loop DC/AC converter in error coordinates with $\mu=\mu_{\textit{ref}}+ d_v(P_{l}-P_{\textit{ref}})$ as%
\begin{align*}
\dot{\tilde \xi} =&~\tilde v_{dc} 
\\
(C_{dc}+&{K_d}) \dot{\tilde v}_{dc}=\!-(G_{dc}\!+\!K_p)\tilde v_{dc} -K_i \tilde \xi-\frac{\mu_{\textit{ref}}}{2}\mathbf{e}_2^{\top}\tilde i_{dq}
\\
&-\!\frac{d_v (P_{l}\!-\!P_{\textit{ref}})}{2}\mathbf{e}_2^{\top} \tilde i_{dq} - \frac{d_v \tilde P_{l}}{2}\mathbf{e}_2^{\top} i_{dq}^{*}
\\
L \dot{\tilde i}_{dq} =&-\!{\left(\mathcal{Z}_{} + \tilde{v}_{dc}\eta L \boldsymbol{J}\right)}  \tilde i_{dq} -\tilde {v}_{dc}\eta L { \boldsymbol{J}}   i_{dq}^* +\frac{\mu_{\textit{ref}}}{2} \mathbf{e}_2  \tilde v_{dc} \nonumber 
\\
&+\frac{d_v(P_{l}\!-\!P_{\textit{ref}})}{2} \mathbf{e}_2  \tilde v_{dc} + \frac{d_v \tilde P_{l}}{2} \mathbf{e}_2  v_{dc}^{*} - \tilde{v}_{dq}
\\
C\, \dot{\tilde{v}}_{dq} =&-\!{\left(\mathcal{Y}_{} + \tilde{v}_{dc}\eta C \boldsymbol{J}\right)}\tilde {v}_{dq}-\tilde{v}_{dc}\eta C {\boldsymbol{J}} v_{dq}^* +\tilde{i}_{dq}
\,,
\end{align*}
with $\tilde P_l=P_l-P_l^*$ and $P_l^*$ as the value of the load power at steady-state.
The derivative of $\mathcal{ {V}}_2$ can be obtained analogously to the proof of Theorem \ref{thm: governor-ctrl},  as
$$\mathcal{\dot V}_2
= -\begin{bmatrix} \tilde v_{dc}& \tilde i_{dq}^\top & \tilde v_{dq}^\top
\end{bmatrix} (\mathcal{Q}  + d_v \mathcal{M})\begin{bmatrix} \tilde v_{dc}& \tilde i_{dq}^\top & \tilde v_{dq}^\top
\end{bmatrix}^\top \leq 0,$$
where $\mathcal{Q}$ is as in \eqref{eq: P-matrix} and $\mathcal{M}$ is a constant matrix, {i.e., its entries do not depend on the droop coefficient $d_v$.} 
Since $\mathcal{Q}$ is positive definite under condition \eqref{eq: neg-def-cond-PID}, there exists $d_v>0$ sufficiently small such that $\mathcal{Q} + d_v \mathcal{M}$ is positive definite. 
A LaSalle-type argument accounting for the evolution of $\tilde \xi$ then asserts global asymptotic stability of the load-induced equilibrium.
\oprocend
\end{pf}

\section{Numerical case study} 
\label{sec: case study}

We validate and test the proposed controllers in a numerical case study. We consider an inverter designed for $10^4$ W power output with the following parameters\footnote{all units are in S.I.}: $G_{dc}=0.1, C_{dc}=0.001,\, R=0.1,\, L=5\cdot 10^{-4},\, C=10^{-5}$, and nominal 
DC voltage of $v_{dc, \textit{ref}} = v_{dc}(0)= 1000$. In order to obtain the desired open-circuit (no load) values $\norm{v_{x}^*}= r_\textit{ref} =165$ and $\omega^*=\omega_0=2\pi50$, we choose the constant gains 
$\eta = \tfrac{\omega_0}{v_{dc,\textit{ref}}}=0.3142,\, 
\mu =\tfrac{2 r_{\textit{ref}}}{v_{dc,\textit{ref}}}=0.33.$

\subsection{Voltage and frequency regulation -- single inverter}
To validate our results for frequency and amplitude regulation, we implement the matching control \eqref{eq: matching control} and the frequency regulation \eqref{eq: idc-ctrl}, together with the three different amplitude controllers. 
We consider a load step of 55\% at $t=0.5s$. The resulting amplitudes and power waveforms are shown in Figure \ref{fig: single-case-zoomed}, whereas Figure \ref{fig: filter-output} shows a time-domain electromagnetic transient (EMT) simulation of the output capacitor voltage.  

The parameters of the frequency controller \eqref{eq: idc-ctrl} were selected as $i_{dc,\textit{ref}}=100, K_p=1, K_i=10,\, K_d=0$ and $\xi(0)=0$. For voltage control we consider the feedforward control \eqref{eq:mu-root}, PI-PBC \eqref{eq: amp-PBCc-ctrl} with (in S.I $\kappa_p=0.1,\kappa_i=10,\, \nu(0)=0$), as well as droop control \eqref{eq: amp-droop-ctrl} (in S.I $\mu_{\textit{ref}}=0.33,\, d_v=10^{-5}$ and $P_{\textit{ref}}=10^4$) plotted as red, green and blue signals, respectively. For all considered controllers, the DC voltage exactly tracks the reference voltage $v_{dc,\textit{ref}}=1000$. The feedforward and PI-PBC designs also track the desired amplitude $r_{\textit{ref}}=165$. Observe that the constant amplitude objective of these controllers requires higher steady-state current amplitudes after the load step. The droop controller on the other hand ensures a trade-off between the power load and AC voltage amplitude. We observe that all controllers yield well-behaved transient response to the step in disturbance.

\subsection{Multi-Converter Case Study}
Next we consider a network of two inverters connected in parallel to a conductance load via a $\Pi$-transmission line model; see Figure \ref{fig: star-topo}. 

\begin{figure}[!ht] 
\centering
\resizebox{8.5cm}{2.5cm}{
\begin{circuitikz}[american voltages]
\draw
(0,0) to [short, *-] (3,0)
(6.5,0) to [R, l=$G_{net}$] (6.5,3) 
(4.7,3) to [open, v^=${v}_{net}$] (4.7,0)
   (6.5,3) to [short, *-, i=$ $] (6.5 ,2)
(7.4,3) to [short, *-, i=$ $] (7.4 ,2)
(7.5,3) to [open, v^=$ $] (7.5,0)
(7.4,3) to [C, l=$C_{net}$] (7.4,0)
(6.7,3) -- (7.7,3) 
(4,0) to [short] (5,0)
(3,0) to (4,0)
(3,3) to (4,3)
(0.5,3) to [open, v^=${v}_{x,1}$] (0.5,0) 
(0,3) 
to [short,*-, i=$i_{\alpha\beta,1}$] (1,3) 
to [R, l=$R$] (2,3) 
to [L, l=$L$] (3.3,3) 
to [short,*-, i_=$ $] (3.3,2) 
to [C, l=$ $] (3.3,1) 
to [short] (3.3,0)
(2.0,3) to [open, v^=$ v_{\alpha\beta,1}$] (2.0,0); 
\draw
(4.1,3) to [short, i<=$i_{net,1}$] (4.0,3) 
(4,3)to [R, l=$R_{net}$] (5.2,3) 
(5.2,3)to [L, l=$L_{net}$] (5.8,3) 
to [short] (7,3);
\draw
(9.7,3) to [short, i=$i_{net,2}$] (9.6,3) 
(7.7,3)to [R, l=$R_{net}$] (8.5,3) 
to [L, l=$L_{net}$] (9.5,3) 
to [short] (10.5,3);
\draw 
(10.5,3) to [short,*-, i_=$ $] (10.5,2) 
(10.5,2)to [C, l=$ $] (10.5,1) 
to [short] (10.5,0);
\draw
(8.3,0)	to [short] (10.55,0)
(10.7,3) to [open, v^=$ v_{\alpha\beta,2}$] (10.7,0); %
\draw
(8.3,0) to [short] (5,0)
(13.5,3) to [short,*-, i>=$i_{\alpha\beta,2}$] (12.5,3) 
(10.5,3)to [R, l=$R$] (11.5,3) 
(11.5,3)to [L, l=$L_{}$] (12.5,3) 
(12.5,3) to [open, v^=${v}_{x,2}$] (12.5,0)
(13.5,0) to [short, *-] (10.5,0);
\end{circuitikz}
}
\caption{Two inverters connected in parallel to a conductance load  $G_{net}>0$ via a $\Pi$-line model.} 
\label{fig: star-topo}
\end{figure}
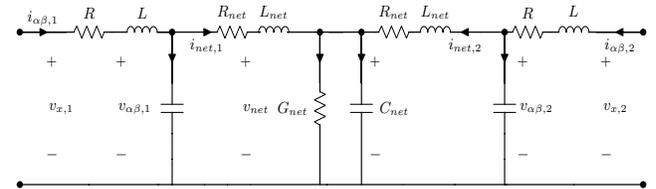

\begin{figure*}[t]
\centering
\includegraphics*[scale=0.5]{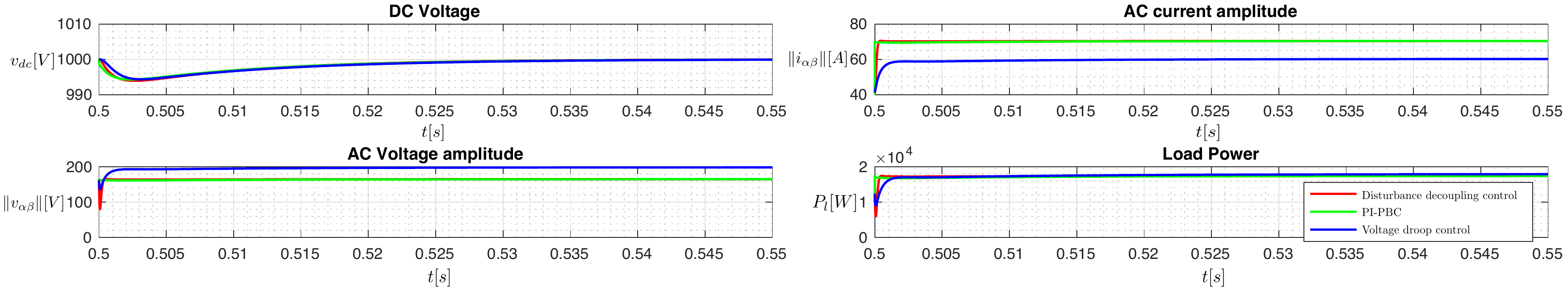}
\caption{The response of the system under the three controllers from Section \ref{sec: outer loops} after a step-up in load conductance at $t=0.5$.}
\label{fig: single-case-zoomed}
\end{figure*}

\begin{figure*}[t]
\centering
\includegraphics*[scale=0.507]{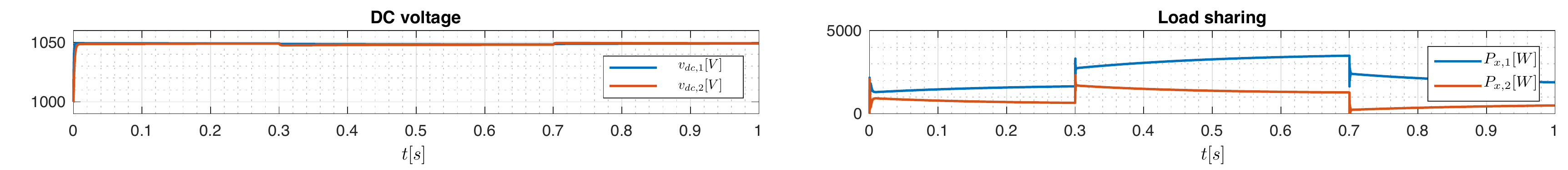}
\caption{The response of the parallel converter scenario in Figure \ref{fig: star-topo} during two steps in the load conductance $G_{net}$.} 
\label{fig: star-topo-sim}
\end{figure*}

\begin{figure}[!h]
\centering
\includegraphics*[width=0.9\columnwidth]{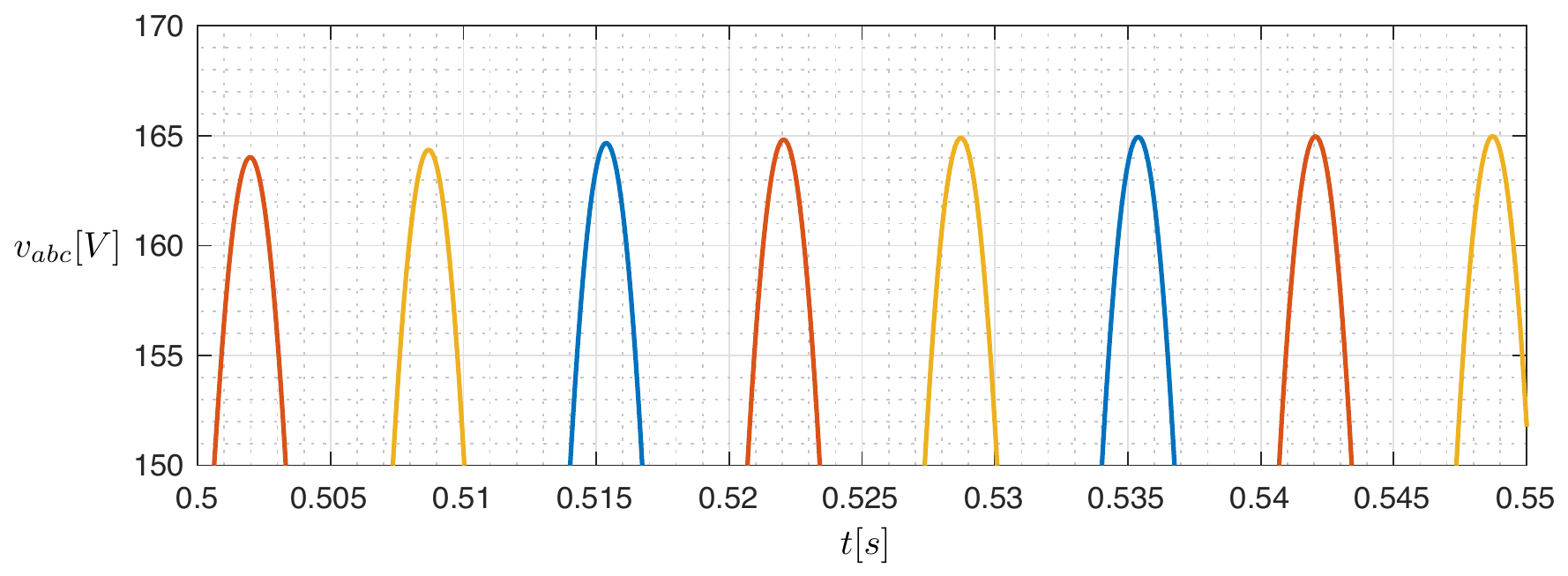}
\caption{Magnified plot of the three-phase AC bus voltage.}
\label{fig: filter-output}
\end{figure}

The $\Pi$-line parameters are $R_{net}=0.5 $, $L_{net}=2.5 \cdot 10^{-5}$, and $C_{net}=2\cdot 10^{-7}$, where the capacitors account for filter and line charge capacitance. The interconnection dynamics are considered for $k\in\{1,2\}$:  
\begin{align*}
C \dot v_{\alpha\beta,k}&= -G v_{\alpha\beta,k} + i_{\alpha\beta,k}-i_{net,k}\\
L_{net} \dot {i}_{net,k} &=-R_{net}i_{net,k} +v_{\alpha\beta,k}-v_{load}\\
C_{net} \dot v_{net}&=-G_{net} v_{net}+ i_{net,1}+i_{net,2} \,.
\end{align*}
We implemented the matching control \eqref{eq: matching control} with gains according to \eqref{eq: pow-sharing} to demonstrate the proportional power sharing with ratio $\rho=3$. We chose the current control parameters for the inverters as {$i_{dc,\textit{ref},1}=100, K_{p1}=2$}, neglected internal losses $G_{dc,1}=G_{dc,2}=0$, fixed the modulation amplitude at $\mu_{1}=\mu_{2}=0.33$, removed the integral action $K_{i,1}=K_{i,2}=0$, and set all other parameters as before. Our simulation in Figure \ref{fig: star-topo-sim} displays a prescribed power sharing ratio of 3:1 under resistive load steps at times $t=0.3$ and $t=0.7$.
%

\section{Conclusions}
\label{sec: conclusions}

This paper addresses the problem of designing grid-forming converter control strategies for weak-grid scenarios, those in which no other unit is able to regulate the AC grid frequency.
Based on the idea of matching the dynamics of a SM, we enable by feedback the crucial coupling between the inverter's DC-side voltage and its AC-side frequency. 
%
%
As a result, the AC grid frequency measurement is replaced by that of the DC-link voltage, further obviating the conventional time scale-separation approach.
The seamless compatibility with synchronous machines, set up by the matching control, yields droop and proportional power sharing characteristics, while preserving passivity properties for the inverter. 
Moreover, the addition of synthetic damping and inertia is straightforward.
By paring the proposed controller with additional outer loops, we also study output voltage regulation in the presence of measurable disturbances. 
These outer controllers are designed based on passivity-based and disturbance decoupling methods and achieve exact tracking for two quantities of interest: output voltage frequency and its amplitude. 
In the light of our analysis, a natural counterpart is to investigate the compatibility of networked objectives and to design suitable controllers that encompass multiple converters.

\bibliographystyle{elsarticle-harv} 
\bibliography{bib/aut}   
\appendix

\end{document}